\documentclass[letterpaper,12pt]{article}
\usepackage{amsmath}
\usepackage{amssymb}
\usepackage{geometry}
\usepackage[colorlinks=true,linkcolor=red,citecolor=blue]{hyperref}
 \usepackage{amssymb}
 \usepackage{amsmath}
 \usepackage[latin1]{inputenc}
 \usepackage[T1]{fontenc}
\usepackage{times}
\usepackage{epsfig}
 \usepackage{dsfont}
 \usepackage{natbib}
 \usepackage{mathrsfs}

\newtheorem{theorem}{Theorem}[section]

\newtheorem{corollary}[theorem]{Corollary}

\newtheorem{remark}[theorem]{Remark}

\numberwithin{equation}{section} 

\begin{document}

\title{Uniform-in-bandwidth consistency for kernel-type estimators of Shannon's entropy}

\author{Salim BOUZEBDA
\footnote{e-mail: salim.bouzebda@upmc.fr} \hbox{ } and \hbox{ }Issam
ELHATTAB\footnote{e-mail:
issam.elhattab@upmc.fr}\\
Laboratoire de Statistique Th\'{e}orique et Appliqu\'{e}e (L.S.T.A.)\\
Universit\'{e} de Paris VI\\4 place Jussieu
    75252 Paris Cedex 05, France}
\date{}
\maketitle

\begin{abstract}
We establish uniform-in-bandwidth consistency for kernel-type
estimators of the differential entropy. We consider two kernel-type
estimators of Shannon's entropy. As a consequence,  an asymptotic
$100\%$ confidence interval of entropy is provided.
\\
\noindent\textit{AMS Subject Classification:}  62F12 ; 62F03 ; 62G30 ; 60F17 ; 62E20.
\\
\noindent{\small {\it Keywords:} Entropy; Kernel
 estimation; Uniform in bandwidth; Consistency.}
\end{abstract}

\section{Introduction and estimation}
Let $\left({\bf X}_n\right)_{n\geq1}$ be a
sequence of independent and identically distributed  $\mathbb{R}^d$-valued random vectors, $d\ge1$,
with cumulative distribution function $\mathbb{F}({\bf
x})=\mathbb{P}({\bf X}\leq {\bf x})$ for ${\bf x}\in\mathbb{R}^d$ and a density function $f(\cdot)$
with respect to Lebesgue measure on $\mathbb{R}^d$.
Here, as usual, ${\bf X}=(X_1,\dots,X_d)\leq {\bf
x}=(x_1,\dots,x_d)$ means that each component of
${\bf X}$ is less  than or equal to the corresponding component of
${\bf x}$, that is, $X_i\leq x_i$, for all $i=1,\dots,d$. The differential (or Shannon)
entropy of $f(\cdot)$ is defined to be
%e2 ###
%e1 ###
\begin{eqnarray}\label{Chap3-entropy}
H(f) &:=& -\int_{\mathbb{R}^d}f({\bf x})\log\left(f({\bf
x})\right)d{\bf x}\\
&:=&-\int_{\mathbb{R}^d}\log\left(f({\bf
x})\right)d\mathbb{F}({\bf
x})\label{secon},
\end{eqnarray}
whenever this integral is meaningful, and where, for $\mathbf{ x}=(x_{1},\ldots,x_{d})$,  $d{\bf x}$ denotes Lebesgue measure in $\mathbb{R}^d$.
We will use the convention that $0\log(0)=0$ since $u\log(u)\rightarrow0$ as $u\rightarrow0$.

 The concept of
differential entropy was originally introduced in Shannon's paper
\cite{shannon1948}. Since this early epoch, the notion of
entropy has been the subject of great theoretical and applied
interest. We refer to \cite[Chapter
8.]{coverthomas2006}  for a comprehensive overview of differential entropy and their
mathematical properties. Entropy concepts and principles play an
fundamental  role in many applications, such as statistical
communication theory \cite{Gallager1968}, quantization theory
\cite{Reny1959}, statistical decision theory
\cite{kullback1959}, and  contingency table analysis
\cite{Gokhale1978}. \cite{Csizar1962} introduced the concept of
convergence in entropy and showed that the latter convergence
concept implies convergence in $\mathcal{L}_1$. This property
indicates that entropy is a useful concept to measure
``\emph{closeness in distribution}'', and also justifies
heuristically the usage of sample entropy as test statistics when
designing entropy-based tests of goodness-of-fit. This line of
research has been pursued by
\cite{vasicek1976,Prescott1976,DudewiczMeulen1981,Gokhale1983,Ebrahimi1992}
and \cite{Esteban2001} [including the
references therein].  The idea here is that many families of
distributions are characterized by maximization of entropy subject
to constraints (see, e.g., \cite{Jaynes1957} and
\cite{Verdugo1978}). There is a huge literature on the Shannon's entropy and its applications. It is not the purpose of this paper to survey this extensive literature.

In the literature, various estimator for $H(f)$, based on a random sample $\mathbf{X}_{1},\ldots,$ $\mathbf{X}_{n}$ from the underlying distribution, have been proposed and their asymptotic properties studied. For an exhaustive list of references in this vein,  we refer to
\cite{gyorfimeulen1990,beirlantdudewiczgyorfivandermeulen1997} and the references therein.

We mention that there exist mainly two  approaches to the construction of entropy estimators. The first approach is based on spacings when $d=1$.
The the second approach, to be used in this paper to estimate $H(f)$, consists in first obtaining a suitable density estimate $f_{n}(\cdot)$ for $f(\cdot)$, and then substituting $f(\cdot)$ by $f_{n}(\cdot)$ in an entropy-like functional of $f(\cdot)$.

  The main contribution of the present paper
is to establish an almost sure uniform in bandwidth consistency of the kernel-type estimator of the entropy
functional $H(f)$. In the entropy framework, the results obtained here  are believed to be novel. %\vskip5pt \noindent

We start by giving some notation and conditions that are needed for the forthcoming sections. Below, we will work under the following assumptions on $f(\cdot)$ to establish our results.
\begin{description}
    \item[{\bf(F.1)}] The functional $H(f)$ is well-defined by (\ref{Chap3-entropy}), in the sense that
%e3 ###
            \begin{eqnarray}\label{Chap3-entropyfinite}
                \vert H(f)\vert<\infty.
            \end{eqnarray}
\end{description}
We recall from  (cf. \cite[p. 237]{Ash1965}, \cite[p. 108]{Berger1971})  that the finiteness of $H(f)$ is guaranteed if both $\mathbb{E}\|\mathbf{ X}\|^{2}<\infty$, where $\|\cdot\|$ denotes the Euclidian norm in $\mathbb{R}^{d}$, (in which case $H(f)<\infty$)  and $f(\cdot)$ is bounded (in which case $H(f)>-\infty$). Ash gives an example of a density function  on $\mathbb{R}$ for which $H(f)=\infty$ and also one for which $H(f)=-\infty$. We refer to \cite[Section 4]{gyorfimeulen1991} for conditions characterizing
(\ref{Chap3-entropyfinite}) in terms of $f(\cdot)$.

  To define our
entropy estimator we define, in a first step, a kernel density
estimator. Towards this aim, we introduce a measurable function
$K(\cdot)$ fulfilling the following conditions.
\begin{description}
    \item[(K.1)]  $K(\cdot)$ is of bounded variation
        on $\mathbb{R}^d$;
    \item[(K.2)] $K(\cdot)$ is right continuous on
        $\mathbb{R}^d$, i.e., for any ${\bf t}=(t_1,\dots,t_d)$, we have
      \begin{eqnarray*}
 K(t_1,\dots,t_d)=\lim_{\varepsilon_1\downarrow0,\dots,
        \varepsilon_d\downarrow0}K(t_1+\varepsilon_1,\dots,t_d+\varepsilon_d);  \end{eqnarray*}
    \item[(K.3)] $\|K\|_\infty := \sup_{{\bf t}\in
        \mathbb{R}^d}\vert K({\bf t})\vert=:\kappa< \infty;$
    \item[(K.4)] $\int_{\mathbb{R}^d} K({\bf
        t})d{\bf t}=1$.
\end{description}
The well known Akaike-Parzen-Rosenblatt (refer to
\cite{akaike1954,Parzen1962} and \cite{rosenblatt1956})
kernel estimator of $f(\cdot)$ is defined, for any ${\bf x}\in\mathbb{R}^d$,
by
%e4 ###
\begin{eqnarray}\label{Chap3-A}
  f_{n,h_n}({{\bf x}}):=
 (nh_n^{d})^{-1}\sum_{i=1}^{n} K(({{\bf x}}-{\bf
 X}_i)/h_n),
\end{eqnarray}
where $0<h_n\leq1$ is the smoothing parameter.
For notational convenience, we have chosen the same bandwidth sequence for each margins. This assumption can be dropped easily. Refer for example to \cite[Remark 8]{Mason2005} for more details.

In a second step, given
$f_{n,h_n}(\cdot)$, we
estimate $H(f)$ using the representation (\ref{Chap3-entropy}),  by setting
%e5 ###
\begin{eqnarray}\label{Chap3-integralestimator}
H_{n,h_n}^{(1)}(f) := -\int_{A_{n}}
  f_{n,h_n}({\bf x})\log\big(  f_{n,h_n}({\bf x})\big)
d{\bf x},
\end{eqnarray}
where $$A_{n}:=\{{\bf x} :  f_{n,h_n}({\bf x})\ge
\gamma_{n}\},$$ and $\gamma_{n}\downarrow 0$ is a sequence of positive constant.
The \emph{plug-in} estimator $H_{n,h_n}^{(1)}(f)$  was introduced by \cite{dmitrievtarasenko1973} for $d=1$ and $A_n = [-b_{n},b_{n}]$, where $b_{n}$ is a specified sequence of constants. The integral estimator $H_{n,h_n}^{(1)}(f) $ can  be easily calculated if, for example, $f_{n}(\cdot)$ is a histogram.

 In the present paper, we will consider also the \emph{resubstitution} estimate
proposed in \cite{ahmadlin1976}.
In this case, we shall study uniform-in-bandwidth consistency of the estimator of $H(f)$ based on the representation (\ref{secon}) which is, in turn,  defined by
%e6 ###
\begin{equation}\label{Chap-5-splittingdataestimate}
H_{n,h_n}^{(2)}(f) := -\frac{1}{n}\sum_{i=1}^{n}
\mathds{1}_{\Omega_{n,i}}\log\left(f_{n,h_n}({\bf X}_i)\right),
\end{equation}
where $$\Omega_{n,i}:=\{  f_{n,h_n}({\bf X}_{i})\ge
\gamma_{n}\}, ~~\mbox{for} ~~i=1,\ldots,n$$
The limiting behavior of
$  f_{n,h_n}(\cdot)$, for appropriate choices of the
bandwidth $h_n$,  has been studied by a large number statisticians over many decades. For good sources of references to research literature in this area along with statistical applications consult
 \cite{devroyelugosi2001,DevroyeGyorfi1985,bosqlecoutre1987,scott1992} and \cite{prakasa1983}. In
particular, under our assumptions, the condition that
$h_n\rightarrow 0$ together with $nh_n\rightarrow\infty$ is
necessary and sufficient for the convergence in probability of
$  f_{n,h_n}({\bf x})$ towards the limit $f({\bf x})$,
independently of ${\bf x}\in\mathbb{R}^d$ and the density
$f(\cdot)$. Various uniform consistency results involving the
estimator $  f_{n,h_n}(\cdot)$ have been recently established.
We refer to \cite{Deheuvels2000,Masoneinmahl2000,Deheuvelsmason2004} and the references therein.
The first seminal paper that devoted  to obtaining uniform in bandwidth results for the kernel-type estimator  was  \cite{Mason2005}. Since then, there is a considerable interest
in obtaining so-called uniform in bandwidth results for kernel-type estimators depending on a bandwidth sequence.
%
%They used sophisticated mathematical techniques from other fields such as topology basing on papers by M. Talagrand, and considered aspects like measurability, precisely, the major idea of their proof combines an exponential inequality of \cite{Talagrand1994} in connection with a suitable moment inequality.
%
In this paper we will use their methods to establish convergence results for
the estimates $H_{n,h_n}^{(1)}(f)$ and  $H_{n,h_n}^{(2)}(f)$ of $H(f)$ in the same spirit of \cite{Bouzebda-Elhattab2009,Bouzebda-Elhattab2010}.

 The remainder of this paper is organized as follows. In Section
\ref{section2}, we state our main results concerning the limiting
behavior of $H_{n,h_n}^{(1)}(f)$ and $H_{n,h_n}^{(2)}(f)$. Some
concluding remarks and possible future developments are mentioned
in Section \ref{conclusion}.  To avoid interrupting the flow of the
presentation, all mathematical developments are relegated to Section \ref{proof}.

%s2 ###
\section{Main results}\label{section2}
To prove the strong consistency of
$H_{n,h_n}^{(1)}$, we shall consider another, but more
appropriate and more computationally convenient, centering factor than the expectation
$\mathbb{E}H_{n,h_n}^{(1)}$, which is delicate to handle.
This is given by
\begin{eqnarray*}
\widehat{\mathbb{E}}H_{n,h_n}^{(1)}(f) := -\int_{A_{n}}
\mathbb{E}  f_{n,h_n}({\bf x})\log\big(\mathbb{E}  f_{n,h_n}({\bf x})\big)
d{\bf x}.
\end{eqnarray*}
\noindent The main result, concerning $H_{n,h}^{(1)}$, to be proved here may now be stated precisely as
follows.

\begin{theorem}\label{Chap3-theorem1}
Let $K(\cdot)$ satisfy {\rm (K.1-2-3-4)}, and let $f(\cdot)$ be a bounded density
fulfill {\rm (F.1)}.
Let $c>0$ and  $\{h_n\}_{n\geq 1}$ be a sequence of positive constants such that, $cn^{-1}\gamma_{n}^{-4}(\log n)\leq
h_{n}<1$.  Then there exists a positive
constant $\Upsilon$, such that
\begin{eqnarray*}
\limsup_{n\rightarrow\infty} \sup_{h_{n}\leq
h\leq 1}\frac{\sqrt{nh\gamma_{n}^{4}}\vert H_{n,h}^{(1)}(f) -
\widehat{\mathbb{E}}H_{n,h}^{(1)}(f)\vert}{\sqrt{\left(\log(1/h)\vee\log\log n\right)}}\leq
\Upsilon~~a.s.
\end{eqnarray*}
\end{theorem}
The proof of Theorem \ref{Chap3-theorem1} is postponed  until \S \ref{proof}.

%\vskip7pt\noindent
Let $(h^\prime_n)_{n\geq1}$ and $(h^{\prime\prime}_n)_{n\geq1}$ be
two sequences of constants such that $0<h^\prime_n<h^{\prime\prime}_n<1$, together with
$h^{\prime\prime}_n\rightarrow0$ and  $nh^\prime_n\gamma_{n}^{4}/\log n\rightarrow\infty,$ as $n\rightarrow\infty$.
A direct application of Theorem \ref{Chap3-theorem1} shows that, with probability 1,
\begin{eqnarray*}
\sup_{h^\prime_n\leq h\leq h^{\prime\prime}_n}\vert H_{n,h}^{(1)}(f) -
\widehat{\mathbb{E}}H_{n,h}^{(1)}(f)\vert=O\left(\sqrt{\frac{(\log (1/h^\prime_n)\vee \log \log n)}{nh^\prime_n\gamma_{n}^{4}}}\right).
\end{eqnarray*}
This, in turn, implies that
%e7 ###
\begin{eqnarray}\label{Chap3-lim}
\lim_{n\rightarrow\infty}\sup_{h^\prime_n\leq h\leq h^{\prime\prime}_n}\vert
H_{n,h}^{(1)}(f) - \widehat{\mathbb{E}}H_{n,h}^{(1)}(f)\vert=0~~~a.s.
\end{eqnarray}
The following result handles the uniform deviation of the estimate
$H_{n,h_n}^{(1)}(f)$ with respect to $H(f)$.
\begin{corollary}\label{Chap3-corol1}
Let $K(\cdot)$ satisfy {\rm (K.1-2-3-4)}, and let $f(\cdot)$ be a
uniformly Lipschitz continuous and  bounded density on
$\mathbb{R}^d$, fulfilling {\rm (F.1)}. Then for each pair of sequences
$0<h^\prime_n<h^{\prime\prime}_n\leq 1$ with
$h^{\prime\prime}_n\rightarrow0$, $nh^\prime_n\gamma_{n}^{4}/\log
n\rightarrow \infty$ and $|\log(h^{\prime\prime}_n)|/\log\log n \rightarrow \infty$ as $n \rightarrow \infty$, we have
%e8 ###
\begin{eqnarray}\label{Chap3-confi}
\lim_{n\rightarrow\infty} \sup_{h^\prime_n\leq h\leq h^{\prime\prime}_n}\vert H_{n,h}^{(1)}(f) -
H(f)\vert= 0\qquad a.s.
\end{eqnarray}
\end{corollary}
The proof of Corollary \ref{Chap3-corol1} is postponed  until \S \ref{proof}.

\begin{remark}
 We note that the main problem in using entropy estimates
such as (\ref{Chap3-integralestimator}) is to choose properly the smoothing
parameter $h_n$. The uniform in bandwidth consistency result given in
(\ref{Chap3-confi}) shows that any choice of $h$ between
$h^\prime_n$ and $h^{\prime\prime}_n$ ensures the consistency of
$H_{n,h}^{(1)}(f)$. In other word, the fluctuation of the bandwidth in a small interval do not affect the consistency of the nonparametric estimator of $H(f)$.
\end{remark}
Now, we shall establish another result in a similar direction for a class of compactly supported densities. We need the following additional conditions.
\begin{description}
\item[(F.2)]  $f(\cdot)$ has a compact support  say
    $\mathbb{I}$ and is $s$-time continuously differentiable, and there exists a constant
    $0<\mathfrak{M}<\infty$ such that
\begin{eqnarray*}
\sup_{\mathbf{x} \in\mathbb{I}}\left|\frac{\partial ^s f(\mathbf{x})}{\partial x_1^{j_1}\ldots\partial x_d^{j_d}}\right|\leq \mathfrak{M},~~j_1+\cdots+j_d=s.
\end{eqnarray*}
\item[(K.5)] $K(\cdot)$ is of order
    $s$, i.e., for some constant $\mathfrak{S}\neq 0$,
\begin{eqnarray*}
&&\int_{\mathbb{R}^d}t_1^{j_1}\ldots t_d^{j_d}K(\mathbf{t})d\mathbf{t}=0,~~j_1,\ldots,j_d\geq 0,~~j_1+\cdots +j_d=1,\ldots,s-1,\\
&& \int_{\mathbb{R}^d}|t_1^{j_1}\ldots t_d^{j_d}|K(\mathbf{t})d\mathbf{t}=\mathfrak{S},~~j_1,\ldots,j_d\geq 0,~~j_1+\cdots +j_d=s.
\end{eqnarray*}
\end{description}
Under the condition (F.2), the differential entropy of $f(\cdot)$
may be written as follows
\begin{eqnarray*}
H(f)=-\int_{\mathbb{I}}f({\bf x})\log\left(f({\bf
x})\right)d{\bf x}.
\end{eqnarray*}

\begin{theorem}\label{Chap3-corol2}
Let $K(\cdot)$ satisfy {\rm (K.1-2-3-4-5)}, and let $f(\cdot)$
fulfill {\rm (F.1-2)}. Then for each pair
of sequences $0<h^\prime_n<h^{\prime\prime}_n\leq 1$ with
$h^{\prime\prime}_n\rightarrow0$    and $nh^\prime_n/\log
n\rightarrow \infty$  as $n\rightarrow \infty$, we have, for any $\gamma>0$
\begin{eqnarray*}
\limsup_{n\rightarrow\infty} \sup_{h^\prime_n\leq
h\leq h^{\prime\prime}_n}\frac{\sqrt{nh}\vert H_{n,h}^{(1)}(f) -H(f)\vert}{\sqrt{\log(1/h)\vee\log\log n}}\leq \zeta(\mathbb{I})~~a.s.,
\end{eqnarray*}
where
\begin{eqnarray*}
\zeta(\mathbb{I}):=\left(\frac{\gamma^{2}+\gamma+1}{\gamma^{2}}\right)^{1/2}\sup_{{\bf x}\in \mathbb{I}}\left\{f({\bf x})\int_{\mathbb{R}^d}K^2({\bf u})d{\bf u}\right\}^{1/2}.
\end{eqnarray*}
\end{theorem}
The proof of Theorem \ref{Chap3-corol2} is postponed until \S \ref{proof}.

%\noindent
To state our result concerning $H_{n,h_{n}}^{(2)}(f)$ we need the following additional condition.
\begin{description}
    \item[(F.3)] $\mathbb{E}\Big[\log^2\big(f({\bf
        X})\big)\Big]<\infty$.
\end{description}
\begin{remark}
Condition {\rm (F.3)} is extremely weak and is satisfied by all commonly encountered distributions including many important heavy tailed distributions for which the moments do not exists (see. e.g. \cite{Song2000} for more details and references on the subject.)
\end{remark}
To prove the strong consistency of
$H_{n,h_n}^{(2)}$ we consider the following centering factor
\begin{eqnarray*}
\widehat{\mathbb{E}}H_{n,h_n}^{(2)}(f) := \frac{1}{n}\sum_{i=1}^{n}\mathds{1}_{\Omega_{n,i}}\log\left(\mathbb{E}(f_{n,h_n}(\mathbf{x})\mid{\bf X}_i=\mathbf{x})\right).
\end{eqnarray*}
The main results concerning $H_{n,h}^{(2)}(f)$ are summarized in the following Theorems.
\begin{theorem}\label{Chap3-theo}
Let $K(\cdot)$ satisfy {\rm (K.1-2-3-4)}, and let $f(\cdot)$ be a bounded density
fulfilling {\rm (F.1)}.
Let $c>0$ and  $\{h_n\}_{n\geq 1}$ be a sequence of positive constants such that, $cn^{-1}\gamma_{n}^{-2}(\log n)\leq
h_{n}<1$.  Then there exists a positive
constant $\Upsilon^{\prime}$, such that
\begin{eqnarray*}
\limsup_{n\rightarrow\infty} \sup_{h_{n}\leq
h\leq 1}\frac{\sqrt{nh\gamma_{n}^{2}}\vert H_{n,h}^{(2)}(f) -
\widehat{\mathbb{E}}H_{n,h}^{(2)}(f)\vert}{\sqrt{\left(\log(1/h)\vee\log\log n\right)}}\leq
\Upsilon^{\prime}~~a.s.
\end{eqnarray*}
\end{theorem}
The proof of Theorem \ref{Chap3-theo} is postponed  until \S \ref{proof}.

\begin{theorem}\label{Chap-5-theorem1}
Assume that the kernel function $K(\cdot)$ is compactly supported and satisfies the conditions  {\rm (K.1-2-3-4-5)}. Let $f(\cdot)$ be a bounded density
function fulfilling the conditions {\rm(F.1-2-3)}.
Let $\{h_n^\prime\}_{n\geq 1}$ and $\{h_n^{\prime\prime}\}_{n\geq 1}$
such that $h_n^\prime=An^{-\delta}$ and $h_n^{\prime\prime}=Bn^{-\delta}$
with arbitrary choices of $0<A<B<\infty$ and $(1/(d+4))\leq\delta <1$. Then, for $\gamma>0$,  we have with probability one,
%e9 ###
\begin{equation}\label{Chap-5-result}
\limsup_{n\rightarrow\infty} \sup_{h_n^{\prime}\leq h\leq
h_n^{\prime\prime}} \frac{\sqrt{nh\gamma_{n}^{2}}\vert
H_{n,h}^{(2)}(f)-H(f)\vert}{\sqrt{2\log(1/h)}}\leq\sigma_\mathbb{I},
\end{equation}
where
\begin{equation*}
\sigma_\mathbb{I}:=\frac{1}{\gamma}\left\{\sup_{{\bf x}\in \mathbb{I}}f({\bf x})
\int_{\mathbb{R}^d} K^2({\bf u})d{\bf u}\right\}^{1/2},
\end{equation*}
where $\mathbb{I}$ is given in {\rm (F.2)}.
\end{theorem}
The proof of Theorem \ref{Chap-5-theorem1} is postponed until \S \ref{proof}.

\begin{remark}
Theorem \ref{Chap3-corol2} leads,
using the techniques developed in \cite{Deheuvelsmason2004}, to the construction of
asymptotic $100\%$ certainty interval for the true entropy $H(f)$, i.e., as $n\rightarrow \infty$, for each $\varepsilon>0$
\begin{eqnarray*}
\mathbb{P}\left(H(f)\in\left[H_{n,h}^{(1)}(f)-A_{n,\varepsilon},
H_{n,h}^{(1)}(f)+B_{n,\varepsilon}\right]\right)\approx 100\%,
\end{eqnarray*}
see (\ref{Chap3-c-i-equa3}) bellow for explicit expressions of $A_{n,\varepsilon}$ and $B_{n,\varepsilon}$.
We give
in what follows, the idea how to construct this interval.
Throughout, we let $h\in[h^\prime_n,h^{\prime\prime}_n]$, where
$h^\prime_n$ and $h^{\prime\prime}_n$ are as in Theorem
\ref{Chap3-corol2}. We infer from Theorem \ref{Chap3-corol2}
that, for suitably chosen data-dependent functions
$L_n=L_n(X_1,\dots,X_n)>0$, for each $0<\varepsilon<1$, we have, as $n\rightarrow\infty$,
%e10 ###
\begin{eqnarray}\label{Chap3-c-i-equa}
\mathbb{P}\left(\frac{1}{L_n}\vert H_{n,h}^{(1)}(f)-H(f)\vert
\geq 1+\varepsilon\right)\rightarrow 0.
\end{eqnarray}
\noindent Assuming the validity of the statement (\ref{Chap3-c-i-equa}), we
obtain asymptotic certainty interval for $H(f)$ in the following
sense. For each $0<\varepsilon<1$, we have, as
$n\rightarrow\infty$,
%e11 ###
\begin{eqnarray}\label{Chap3-c-i-equa3}
\mathbb{P}\left(H(f)\in\left[H_{n,h}^{(1)}(f)-(1+\varepsilon)L_n,
H_{n,h}^{(1)}(f)+(1+\varepsilon)L_n\right]\right)\rightarrow 1.
\end{eqnarray}
\noindent Whenever (\ref{Chap3-c-i-equa3}) holds for each
$0<\varepsilon<1$, we will say that the interval
\begin{eqnarray*}
\big[H_{n,h}^{(1)}(f)-L_n,
H_{n,h}^{(1)}(f)+L_n\big],
\end{eqnarray*}
 provides asymptotic $100\%$ certainty
interval for $H(f)$.

 To construct $L_n$ we proceed
as follows. Assume that there exists a sequence $\{\mathbb{I}_n\}_{n\geq1}$ of strictly
nondecreasing compact subsets of $\mathbb{I}$, such that $$\bigcup_{n \geq 1}
\mathbb{I}_n=\mathbb{I}$$ {\rm(}for the estimation of the support
$\mathbb{I}$ we may refer to \cite{DevroyeWise1980} and the references
therein{\rm)}. Furthermore, suppose that there exists a sequence
{\rm(}possibly random{\rm)}
$\{\zeta_n(\mathbb{I}_n)\},n=1,2,\dots,$ converging to
$\zeta(\mathbb{I})$ in the sense that
%e12 ###
\begin{eqnarray}\label{Chap3-unif-equa}
    \mathbb{P}\left(\left\vert
    \frac{\zeta_n(\mathbb{I}_n)}{\zeta(\mathbb{I})}-1\right\vert\geq
    \varepsilon\right)\rightarrow 0 \qquad\mbox{as}~~n\rightarrow\infty~~
    \mbox{for each}~~\varepsilon>0.
\end{eqnarray}
Observe that the statement (\ref{Chap3-unif-equa}) is satisfied when the choice
$$\zeta_n(\mathbb{I}_n):=\sup_{{\bf x}\in
\mathbb{I}_n} \sqrt{  f_{n,h}({\bf
x})\int_{\mathbb{R}^d}K^2({\bf u})d{\bf u}}$$ is considered. Consequently,
we may define the quantity $L_n$ displayed in the statement (\ref{Chap3-c-i-equa}) by
\begin{eqnarray*}L_n:=\sqrt{\frac{\gamma_{n}^{4}\big(\log(1/h)
\vee\log\log n\big)}{nh}}\times\zeta_n(\mathbb{I}_n).
\end{eqnarray*}

\end{remark}
\begin{remark} A practical choice of $\gamma_{n}$ is $\beta (\log n)^{-\alpha}$ where  $\beta>0$ and $\alpha\geq 0$. In the case of the density which is bounded away from $0$, $\alpha$ is equal to $0$.
\end{remark}

\begin{remark}
  \cite{mason2008} establish uniform in bandwidth consistency and central limit theorems for a different but related estimator
to the one proposed in the present paper. That is, \cite{mason2008} propose
  \begin{eqnarray*}
        \widehat{H}_{n,h_n} := -\frac{1}{n}\sum_{i=1}^n\log\left\{f_{n,h_n,-i}(X_i)\right\},
  \end{eqnarray*}
  where
  \begin{eqnarray*}
f_{n,h_n,-i}(X_i):=1/((n-1)h_n)\sum_{1\leq j\neq i\leq
n}K\left((X_i-X_j)/h_n\right).
  \end{eqnarray*}
Their results hold subject to the condition ((C) p. 751, where we choose $\phi(x)=x\log x$ that corresponds to the negative entropy) which is satisfied when density $f(\cdot)$ is
   bounded away from $0$ on its support, refer to Remark 1. p. 752 of \cite{mason2008}, their approach is different from that used in this paper and is based on the notion of a local $U$-statistic. We mention that the estimator proposed by \cite{mason2008} seems to be simpler and with faster rates of convergence.
   The fact that we use the a ``thresholding'' estimator of the entropy permits us  to consider a large class of density by paying the price of loss in  the rate of convergence. Furthermore, if we assume that the density $f(\cdot)$ is bounded
away from $0$ on its support, then the rate of the strong
convergence is of order $\{\{\log(1/h_n)\}/\{nh_n\}\}^{1/2}$ which
is the same rate of the strong convergence for the density
kernel-type estimators, this is precisely the contain of  Theorem \ref{Chap3-corol2}.
\end{remark}

%s3 ###
\section{Concluding remarks and future works}\label{conclusion}
We have addressed the problem of nonparametric estimation of Shannon's entropy.  The results presented in this work are general, since the
required conditions are fulfilled by a large class of densities.

% \vskip5pt \noindent
The evaluation of the integral in (\ref{Chap3-integralestimator}) requires numerical integration and is not easy if $f_{n,h_{n}}(\cdot)$ is a kernel density estimator but it does not involve any
stochastic aspects. The integral estimator can however be easily calculated if we approximate $f_{n,h_{n}}(\cdot)$  by piecewise-constant functions on a
fine enough partition, for example, $f_{n,h_{n}}(\cdot)$ is a histogram.
We mention that in some particular case ($K(\cdot)$ is a double exponential kernel), the
approximations are easily calculated since the distribution function
corresponding to the kernel $K(\cdot)$ is available, confer \cite{EggermontLariccia1999} for more details.
An interesting aspect of the $H_{n,h_n}^{(2)}(f)$ is  that its rate of convergence is faster than that of $H_{n,h_n}^{(1)}(f)$ and that is very easy to compute.

It will be interesting  to enrich our results presented here by an additional uniformity in term of $\gamma_{n}$ in the supremum appearing in all our theorems, which requires non trivial mathematics, this
would go well beyond the scope of the present paper. Another direction of research is to  obtain results, based on $U$-statistic approach, similar to that in \cite{mason2008} for entropy estimator under general conditions, i.e., without assuming the condition that the density $f(\cdot)$ is
   bounded away from $0$ on its support.

%s4 ###
\section{Proofs}\label{proof}
This section is devoted to the proofs of our results.%\vskip5pt

\subsection*{Proof of Theorem \ref{Chap3-theorem1}.} We first
decompose $H_{n,h_n}^{(1)}(f) -
\widehat{\mathbb{E}}H_{n,h_n}^{(1)}(f)$ into the sum of two
components, by writing
%e13 ###
\begin{eqnarray}\label{Chap3-decomposition}
\lefteqn{H_{n,h_n}^{(1)}(f) - \widehat{\mathbb{E}}H_{n,h_n}^{(1)}(f)}\nonumber\\
&=&-\int_{A_{n}}
  f_{n,h_n}({\bf x})\log\big(  f_{n,h_n}({\bf x})\big)
d{\bf x}\nonumber\\
&&+\int_{A_{n}}
\mathbb{E}  f_{n,h_n}({\bf x})\log\big(\mathbb{E}  f_{n,h_n}({\bf x})\big)
d{\bf x}\nonumber\\
%\end{eqnarray}
%\begin{eqnarray}
&=&-\int_{A_{n}} \left\{\log
  f_{n,h_n}({\bf x})-\log\mathbb{E}  f_{n,h_n}({\bf x})\right\}\mathbb{E}  f_{n,h_n}({\bf x})
d{\bf x}\nonumber\\
&&-\int_{A_{n}}
\left\{  f_{n,h_n}({\bf x})-\mathbb{E}  f_{n,h_n}({\bf x})\right\}
\log  f_{n,h_n}({\bf x}) d{\bf x}\nonumber\\
&:=&\boldsymbol{\Delta}_{1,n,h_n}+\boldsymbol{\Delta}_{2,n,h_n}.
\end{eqnarray}
We observe that for all $z>0$, $\left\vert\log z\right\vert \leq
\left\vert\frac{1}{z}-1\right\vert+\left\vert z-1\right\vert$.
Therefore, for any $\mathbf{ x}\in A_{n}=\{{\bf x} :  f_{n,h_n}({\bf
x})\ge \gamma_{n}\}$, we get
 \begin{eqnarray*}
  \lefteqn{\vert\log  f_{n,h_n}({\bf x})-\log \mathbb{E}  f_{n,h_n}({\bf x})\vert =
  \left\vert\log \frac{  f_{n,h_n}({\bf x})}{\mathbb{E}  f_{n,h_n}({\bf x})}\right\vert} \\
   &\leq& \left\vert\frac{\mathbb{E}  f_{n,h_n}({\bf x})}{  f_{n,h_n}({\bf x})}-1\right\vert+
   \left\vert\frac{  f_{n,h_n}({\bf x})}{\mathbb{E}  f_{n,h_n}({\bf x})}-1\right\vert\\
   &=&\frac{\left\vert\mathbb{E}  f_{n,h_n}({\bf x})-
  f_{n,h_n}({\bf x})\right\vert}{  f_{n,h_n}({\bf x})}+
\frac{\left\vert  f_{n,h_n}({\bf x})-\mathbb{E}  f_{n,h_n}({\bf x})\right\vert}{\mathbb{E}  f_{n,h_n}({\bf x})}.
\end{eqnarray*}
In the following $\|\cdot\|_{\infty}$ denotes, as usual, the supremum norm, i.e., $
\|\phi(\mathbf{ x})\|_{\infty}:=\sup_{\mathbf{ x}\in \mathbb{R}^{d}}\|\phi(\mathbf{ x})\|.$
We know (see, e.g., \cite{Mason2005}), for each $h^\prime_n<h<h^{\prime\prime}_n$, as
$n\rightarrow\infty$, we have
\begin{eqnarray*}
\Vert f_{n,h}({\bf x})-\mathbb{E}  f_{n,h}({\bf x})
\Vert_\infty=O\left(\sqrt{\frac{(\log (1/h^\prime_n)\vee \log \log n)}{nh^\prime_n}}\right).
\end{eqnarray*}
For any ${\bf x}\in A_{n}$, one  can see that
%\begin{eqnarray*}
%\mathbb{E}  f_{n,h_{n}}(\mathbf{x})&\geq&  f_{n,h_n}({\bf x})-|  f_{n,h_n}({\bf x})-\mathbb{E}  f_{n,h_{n}}(\mathbf{x})|\geq  \gamma_{n}\\&& +O\left(\sqrt{\frac{(\log (1/h^\prime_n)\vee \log \log n)}{nh^\prime_n}}\right),
%\end{eqnarray*}thus, for $n$ enough large, the second term of the last inequality is dominated by the first one, then, we obtain
\begin{eqnarray*}
\mathbb{E}  f_{n,h_{n}}(\mathbf{x})\geq \gamma_{n}.
\end{eqnarray*}
We readily obtain from these
relations, for any ${\bf x}\in A_{n}$, that
\begin{eqnarray*}
  \vert\log  f_{n,h_n}({\bf x})-\log \mathbb{E}  f_{n,h_n}({\bf x})\vert &\leq&
\frac{2}{\gamma_{n}}\left\vert  f_{n,h_n}({\bf x})-\mathbb{E}  f_{n,h_n}({\bf x})\right\vert.
\end{eqnarray*}
We can therefore write, for any $n\ge1$, the following chain of inequalities
\begin{eqnarray*}
\left\vert\boldsymbol{\Delta}_{1,n,h_n}\right\vert &=&
\left\vert\int_{A_{n}} \left\{\log
  f_{n,h_n}({\bf x})-\log
\mathbb{E}  f_{n,h_n}({\bf x})\right\}\mathbb{E}  f_{n,h_n}({\bf x})
d{\bf x}\right\vert\\
&\leq& \int_{A_{n}}\left\vert\log
  f_{n,h_n}({\bf x})-\log
\mathbb{E}f_{n,h_n}({\bf x})\right\vert
\mathbb{E}f_{n,h_n}({\bf x})d{\bf x}\\
&\leq&\frac{2}{\gamma_{n}}\int_{A_{n}}\left\vert  f_{n,h_n}({\bf x})-\mathbb{E}  f_{n,h_n}({\bf x})\right\vert
\mathbb{E}  f_{n,h_n}({\bf x}) d{\bf x}\\
&\leq&\frac{2}{\gamma_{n}} \sup_{{\bf x}\in A_{n}}
\left\vert \mathbb{E}  f_{n,h_n}({\bf x})-
  f_{n,h_n}({\bf x})\right\vert
\int_{A_{n}} \mathbb{E}  f_{n,h_n}({\bf x})d{\bf x}\\
&\leq& \frac{2}{\gamma_{n}} \sup_{{\bf x}\in
\mathbb{R}^d} \left\vert \mathbb{E}  f_{n,h_n}({\bf x})-
  f_{n,h_n}({\bf x})\right\vert \int_{\mathbb{R}^d}
\mathbb{E}  f_{n,h_n}({\bf x})d{\bf x}.
\end{eqnarray*}
In view of condition (K.4),  by the change of variables and an application of Fubini's theorem, we have
\begin{eqnarray*}
\int_{\mathbb{R}^d}
\mathbb{E}  f_{n,h}({\bf x})d{\bf x}=1.
\end{eqnarray*}
Thus, for any $n\ge1$, we have the following bound
%e14 ###
\begin{eqnarray}\label{Chap3-res-delta1}
\left\vert\boldsymbol{\Delta}_{1,n,h_n}\right\vert \leq \frac{2}{\gamma_{n}} \sup_{{\bf x}\in \mathbb{R}^d} \left\vert
\mathbb{E}  f_{n,h_n}({\bf x})-  f_{n,h_n}({\bf x})
\right\vert.
\end{eqnarray}
We next evaluate the second term $\boldsymbol{\Delta}_{2,n,h_n}$ in the
right side of (\ref{Chap3-decomposition}). Since $\left\vert\log
z\right\vert \leq\frac{1}{z}+z$, for all $z>0$, one can see that
\begin{eqnarray*}
\vert\boldsymbol{\Delta}_{2,n,h_n}\vert &= & \left\vert\int_{A_{n}}
\left\{  f_{n,h_n}({\bf x})-\mathbb{E}  f_{n,h_n}({\bf x})\right\}
\log  f_{n,h_n}({\bf x}) d{\bf x}\right\vert\\
&\leq&\int_{A_{n}}\vert  f_{n,h_n}({\bf x})-
\mathbb{E}  f_{n,h_n}({\bf x})\vert\left[\frac{1}{  f_{n,h_n}({\bf x})}+
  f_{n,h_n}({\bf x})\right]d{\bf x}.
\end{eqnarray*}
Similarly as above, we get, for any ${\bf x}\in A_{n}$,
\begin{eqnarray*}
\frac{1}{  f_{n,h_n}({\bf x})}+
  f_{n,h_n}({\bf x})&=&\left(\frac{1}{  f_{n,h_n}({\bf x})  f_{n,h_n}({\bf x})}+
1\right)  f_{n,h_n}({\bf x})\\
&\leq&\Big(\frac{1}{\gamma_{n}^{2}}+1\Big)  f_{n,h_n}({\bf x}).
\end{eqnarray*}
We can therefore write the following chain of inequalities, for any $n\ge1$,
\begin{eqnarray*}
\lefteqn{\vert\boldsymbol{\Delta}_{2,n,h_n}\vert}\\ &\leq& \Big(\frac{1}{\gamma_{n}^{2}}+1\Big)\int_{A_{n}}\vert
\mathbb{E}  f_{n,h_n}({\bf x})-  f_{n,h_n}({\bf x})\vert\
  f_{n,h_n}({\bf x})d{\bf x}\\
&\leq&\Big(\frac{1}{\gamma_{n}^{2} }+1\Big)\sup_{{\bf x}\in
A_{n}} \left\vert \mathbb{E}  f_{n,h_n}({\bf x})-
  f_{n,h_n}({\bf x})\right\vert\int_{A_{n}}
  f_{n,h_n}({\bf x})d{\bf x}\\
 &\leq&
\Big(\frac{1}{\gamma_{n}^{2} }+1\Big)\sup_{{\bf x}\in
A_{n}} \left\vert \mathbb{E}  f_{n,h_n}({\bf x})-
  f_{n,h_n}({\bf x})\right\vert\int_{\mathbb{R}^d}
  f_{n,h_n}({\bf x})d{\bf x}.
\end{eqnarray*}
In view of condition (K.4), by change of variables, we have
\begin{eqnarray*}
\int_{\mathbb{R}^d}  f_{n,h}({\bf x})d{\bf x}=1.
\end{eqnarray*}
Thus, for any $n\ge1$, we have
%e15 ###
\begin{eqnarray}\label{Chap3-res-delta2}
\vert\boldsymbol{\Delta}_{2,n,h_n}\vert\leq\Big(\frac{1}{\gamma_{n}^{2}}+1\Big)\sup_{{\bf x}\in \mathbb{R}^d}
\left\vert \mathbb{E}  f_{n,h_n}({\bf x})-
  f_{n,h_n}({\bf x})\right\vert.
\end{eqnarray}
%\vskip5pt \noindent
We now impose some slightly more general
assumptions on the kernel $K(\cdot)$ than that of Theorem
\ref{Chap3-theorem1}. Consider the class of functions
\begin{eqnarray*}
\mathcal{K} := \Big\{K(({\bf x}-\cdot)/h^{1/d}): h>0,~{\bf x}\in
\mathbb{R}^d\Big\}.
\end{eqnarray*}
For $\varepsilon>0$, set $N(\varepsilon,\mathcal{K}) = \sup_Q
N(\kappa\varepsilon,\mathcal{K},d_Q)$, where the supremum is taken
over all probability measures $Q$ on $(\mathbb{R}^d,\mathcal{B})$,  where  $\mathcal{B}$ represents the $\sigma$-field of Borel sets of $\mathbb{R}^d$.
Here, $d_Q$ denotes the $L_2(Q)$-metric and
$N(\kappa\varepsilon,\mathcal{K},d_Q)$ is the minimal number of
balls $\{g:d_Q(g,g')<\varepsilon\}$ of $d_Q$-radius $\varepsilon$
needed to cover $\mathcal{K}$. We assume that $\mathcal{K}$
satisfies the following uniform entropy condition.
\begin{description}
    \item[{\bf(K.6)}]   for some $C>0$ and $\nu>0$,
%e16 ###
        \begin{eqnarray}\label{Chap3-entropi}
        N(\varepsilon,\mathcal{K})\leq C\varepsilon^{-\nu},
        0<\varepsilon<1.
        \end{eqnarray}
\end{description}%\vskip5pt
 Finally, to avoid using
        outer probability measures in all of statements, we impose the
        following measurability assumption.
\begin{description}
    \item[{\bf(K.7)}] $\mathcal{K}$ is a pointwise measurable
        class, that is, there exists a countable subclass
        $\mathcal{K}_0$ of $\mathcal{K}$ such that we can find
        for any function $g\in\mathcal{K}$ a sequence
        of functions $\{g_m: m\ge 1\}$ in
        $\mathcal{K}_0$ for which
        \begin{eqnarray*}
          g_m({\bf z})\longrightarrow g({\bf z}),\qquad{\bf z}\in \mathbb{R}^d.
        \end{eqnarray*}
\end{description}
\begin{remark}
Remark that condition (K.6) is satisfied whenever (K.1) holds, i.e., $K(\cdot)$ is of
bounded variation on $\mathbb{R}^d$  (in the sense of Hardy and
Kauser, see, e.g. \cite{Clarkson1933,Vituvkin1955} and
\cite{Hobson1958}). Condition (K.7) is satisfied whenever (K.2) holds, i.e.,
$K(\cdot)$ is right continuous (refer to \cite{Deheuvelsmason2004}
and \cite{Mason2005} and the references therein).
\end{remark}

%\vskip5pt \noindent
By Theorem $1$ of \cite{Mason2005}, whenever
$K(\cdot)$ is measurable and satisfies {\rm (K.3-4-6-7), and when
$f(\cdot)$ is bounded, we have for each $c>0$, and for a suitable
function $\Sigma(c)$, with probability 1,
%e17 ###
\begin{eqnarray}\label{Chap3-resulta-em}
\limsup_{n\rightarrow\infty} \sup_{cn^{-1}\log n\leq
h\leq1}\frac{\sqrt{nh}\Vert
  f_{n,h}-\mathbb{E}  f_{n,h}\Vert_\infty}{\sqrt{\log(1/h)\vee
\log\log n}}=\Sigma(c)<\infty,
\end{eqnarray}
which implies, in view of (\ref{Chap3-res-delta1}) and
(\ref{Chap3-res-delta2}), that, with probability 1,
%e18 ###
\begin{eqnarray}\label{Chap3-result1}
\limsup_{n\rightarrow\infty} \sup_{h_{n}\leq
h< 1}\frac{\sqrt{nh\gamma_{n}^{4}}\vert
\boldsymbol{\Delta}_{1,n,h}\vert}{\sqrt{\left(\log(1/h)\vee\log\log n\right)}}=0,
\end{eqnarray}
and
%e19 ###
\begin{eqnarray}\label{Chap3-result2}
\limsup_{n\rightarrow\infty} \sup_{h_{n}\leq
h< 1}\frac{\sqrt{nh\gamma_{n}^{4}}\vert
\boldsymbol{\Delta}_{2,n,h}\vert}{\sqrt{\left(\log(1/h)\vee\log\log n\right)}}\leq
\Upsilon(c).
\end{eqnarray}
Recalling (\ref{Chap3-decomposition}), the proof of Theorem
\ref{Chap3-theorem1} is completed by combining
(\ref{Chap3-result1}) with (\ref{Chap3-result2}).\hfill$\blacksquare$

\subsection*{Proof of Corollary \ref{Chap3-corol1}.} Recall
$A_{n}=\{{\bf x} :  f_{n,h_n}({\bf x})\ge \gamma_{n}\}$ and let $A^c_{n}$ the complement of
$A_{n}$ in $\mathbb {R}^d$ (i.e., $A^c_{n}=\{{\bf x} :
  f_{n,h_n}({\bf x})<\gamma_{n}\}$). Observe that
\begin{eqnarray*}
|f(\mathbf{x})|&\geq& |  f_{n,h_n}({\bf x})|-|  f_{n,h_n}({\bf x})-f(\mathbf{x})|\geq  \gamma_{n}\\&& +O\left(\sqrt{\frac{(\log (1/h^\prime_n)\vee \log \log n)}{nh^\prime_n}}\right)+O({h^{\prime\prime}_n}^{1/d}).
\end{eqnarray*}
Keep in mind that $|\log(h^{\prime\prime}_n)|/\log\log n \rightarrow \infty$ as $n \rightarrow \infty$, thus, for $n$ enough large, the two last terms of the last inequality are dominated by the first one, then, we obtain
\begin{eqnarray*}
|f(\mathbf{x})|\geq \gamma_{n}.
\end{eqnarray*}
We repeat
the arguments above with the formal change of $H_{n,h_n}^{(1)}(f)$
by $H(f)$. We show that, for any $n\geq 1$,
%e20 ###
\begin{eqnarray}\label{Chap3-bia}
\nonumber \lefteqn{\vert \widehat{\mathbb{E}}H_{n,h_n}^{(1)}(f) -
H(f)\vert}\\\nonumber&\leq&\left\vert\int_{A^c_{n}}f({\bf
x})\log\big(f({\bf x})\big) d{\bf x}\right\vert\\&&+\frac{1}{\gamma_{n}} \sup_{{\bf x}\in \mathbb{R}^d} \left\vert
\mathbb{E}  f_{n,h_n}({\bf x})- f({\bf x})
\right\vert \nonumber\\
&&+\Big(\frac{1}{\gamma_{n}^{2}}+1\Big)\sup_{{\bf x}\in \mathbb{R}^d}
\left\vert \mathbb{E}  f_{n,h_n}({\bf x})-
f({\bf x})\right\vert.
\end{eqnarray}
It is obvious to see that
\begin{eqnarray*}
  \int_{A^c_{n}}f({\bf
x}) d{\bf x} &\leq&  \int_{\frac{1}{2}f({\bf x})\leq \gamma_{n}}f({\bf x}) d{\bf x}+\int_{  f_{n,h}({\bf
x})\leq\gamma_{n}\leq\frac{1}{2}f({\bf x})}f({\bf
x}) d{\bf x}\\
   &\leq&\int_{\frac{1}{2}f({\bf x})\leq \gamma_{n}}f({\bf x}) d{\bf x}+2\int_{\mathbb{R}^d}\vert
  f_{n,h}({\bf x})- f({\bf x})\vert d{\bf x}.
\end{eqnarray*}
Observe that we have
\begin{eqnarray*}
\mathds{1}_{\{\frac{1}{2}f({\bf x})\leq \gamma_{n}\}} f(\mathbf{x})\leq f(\mathbf{x})
\end{eqnarray*}
and $\mathds{1}_{\{\frac{1}{2}f({\bf x})\leq \gamma_{n}\}} f(\mathbf{x})\rightarrow 0$ as $n \rightarrow \infty$, thus  an application of Lebesgue
dominated convergence theorem gives
%e21 ###
\begin{eqnarray}\label{Chap3-eg2}
\lim_{n\rightarrow\infty}\int_{\frac{1}{2}f({\bf x})\leq \gamma_{n}}f({\bf x}) d{\bf x}=0.
\end{eqnarray}
 Keep in mind that the conditions $h_n\rightarrow 0$ together with
$nh_n\rightarrow\infty$ as $n\rightarrow\infty$, ensure that (see
e.g., \cite{DevroyeGyorfi1985})
\begin{eqnarray*}
    \lim_{n\rightarrow\infty}\int_{\mathbb{R}^d}\vert
     f_{n,h_n}({\bf x})- f({\bf x})\vert d{\bf x}=0\qquad a.s.
\end{eqnarray*}
We need the following instrumental fact due to  \cite[Lemma 3.3. p.40]{Devroye1987} and see also \cite[Proof of Theorem 2.2]{Louani2005} which for convenience and easy reference we state here.

\vskip7pt\noindent {\bf Fact.} Let $[h^\prime_n,h^{\prime\prime}_n]$ be a sequence of deterministic interval, where
$nh^\prime_n\rightarrow\infty$ and
$h^{\prime\prime}_n\rightarrow0$, as $n\rightarrow\infty$. For every $\epsilon >0$, then there exist $n_{0}>0$ and $r>0$ such that
$$
\mathbb{P}\left\{\sup_{h^\prime_n\leq h\leq h^{\prime\prime}_n}\int_{\mathbb{R}^d}\vert
     f_{n,h}({\bf x})- f({\bf x})\vert d{\bf x}>\epsilon\right\}\leq \exp\left\{-rn\epsilon^{2}\right\}, ~~n\geq n_{0}.
$$
A routine application of the Borel-Cantelli lemma implies,
for all $h\in[h^\prime_n,h^{\prime\prime}_n]$ such that
$nh^\prime_n\rightarrow\infty$ and
$h^{\prime\prime}_n\rightarrow0$, as $n\rightarrow\infty$, that
%e22 ###
\begin{eqnarray}\label{Chap3-eg1}
    \lim_{n\rightarrow\infty}\sup_{h^\prime_n\leq h\leq h^{\prime\prime}_n}\int_{\mathbb{R}^d}\vert
     f_{n,h}({\bf x})- f({\bf x})\vert d{\bf x}=0\qquad a.s.
\end{eqnarray}
By combining (\ref{Chap3-eg1}) with (\ref{Chap3-eg2}), we obtain
%e23 ###
\begin{eqnarray}\label{Chap3-res1}
\lim_{n\rightarrow\infty}\sup_{h^\prime_n\leq h\leq
h^{\prime\prime}_n}\int_{A^c_{n}}f({\bf x}) d{\bf x}= 0\qquad
a.s.
\end{eqnarray}
Since the entropy $H(f)$ is finite [by condition (F.1)], the measure
\begin{eqnarray*}
\nu(A):=\int_{A}\vert\log\big( f({\bf x})\big)\vert d\mathbb{F}({\bf x}),
\end{eqnarray*}
is absolutely continuous with  respect to the  measure
$\mu(A)=\int_{A}d\mathbb{F}({\bf x})$, which guaranteed that
%e24 ###
\begin{eqnarray}\label{Chap3-pp}
\lim_{n\rightarrow\infty}\sup_{h^\prime_n\leq h\leq
h^{\prime\prime}_n}\int_{A^c_{n}}f({\bf x})\log\big( f({\bf x})\big) d{\bf
x}= 0\qquad a.s.
\end{eqnarray} Recall that we
have for each $h^\prime_n<h<h^{\prime\prime}_n$, as
$n\rightarrow\infty$,
%e25 ###
\begin{eqnarray}
\Vert \mathbb{E}  f_{n,h}({\bf x})- f({\bf x})
\Vert_\infty=O({h^{\prime\prime}_n}^{1/d}).
\end{eqnarray}
Thus, we have
\begin{eqnarray*}
\lim_{n\rightarrow\infty}\sup_{h^\prime_n\leq h\leq h^{\prime\prime}_n}\gamma_{n}^{-2}\Vert \mathbb{E}  f_{n,h}({\bf x})-
f({\bf x})\Vert_\infty=0.
\end{eqnarray*}
This when combined with (\ref{Chap3-bia}), entails that, as
$n\rightarrow\infty$,
%e26 ###
\begin{eqnarray}\label{Chap3-conv}
\sup_{h^\prime_n\leq h\leq h^{\prime\prime}_n}\Vert \widehat{\mathbb{E}}H_{n,h}^{(1)}(f) -
H(f)\Vert\rightarrow 0.
\end{eqnarray}
Using (\ref{Chap3-res1}) and (\ref{Chap3-conv}) in connection with
(\ref{Chap3-lim})
imply the desired conclusion~(\ref{Chap3-confi}).
\hfill$\blacksquare$

\subsection*{Proof of Theorem \ref{Chap3-corol2}.}
Under
conditions (F.2), (K.5) and using Taylor expansion of order $s$ we
get, for ${\bf x}\in \mathbb{I}$,
\begin{eqnarray*}
\vert \mathbb{E}f_{n,h}({\bf x})- f({\bf x})
\vert=\frac{h^{s/d}}{s!}\left\vert\int\sum_{k_1+\cdots+k_d=s}
 t_1^{k_1}\ldots t_d^{k_d}\frac{\partial^s f(\mathbf{x}-h\theta\mathbf{t})}{\partial
 x_1^{k_1}\ldots \partial x_d^{k_d}} K(\mathbf{t})d\mathbf{t}\right\vert,
\end{eqnarray*}
where $\theta=(\theta_1,\ldots,\theta_d)$ and $0<\theta_i<1$, $
i=,1,\ldots,d$. Thus a straightforward application of Lebesgue
dominated convergence theorem gives, for $n $ large enough,
%e27 ###
\begin{eqnarray}\label{equationbiais}
\sup_{{\bf x}\in \mathbb{I}}\vert \mathbb{E}f_{n,h}({\bf x})- f({\bf x})
\vert=O({h_n^{\prime\prime}}^{s/d}).
\end{eqnarray}
Let $\mathbb{J}$ be a nonempty compact subset of the interior of
$\mathbb{I}$ (say $\mathring{\mathbb{I}}$). First, note that we
have from Corollary 3.1.2. p. 62 of \cite{Viallon2006}
%e28 ###
\begin{eqnarray}\label{Chap3-resulta-em-70}
\limsup_{n\rightarrow\infty} \sup_{h_n^\prime\leq h\leq h_n^{\prime\prime}}\sup_{{\bf x}\in \mathbb{J}}\frac{\sqrt{nh}\vert
  f_{n,h}({\bf x})-f({\bf x})\vert}{\sqrt{\log(1/h)\vee
\log\log n}}=\sup_{{\bf x}\in \mathbb{J}}\left\{f({\bf x})\int_{\mathbb{R}^d}K^2({\bf t})d{\bf t}\right\}^{1/2}.
\end{eqnarray}
Set, for all $n\geq 1$,
%e29 ###
\begin{eqnarray}
\pi_n(\mathbb{J})=\left| \int_{\mathbb{J}}
  f_{n,h_n}({\bf x})\log\big(  f_{n,h_n}({\bf x})\big)
d{\bf x}-\int_{\mathbb{J}}
f({\bf x})\log\big(f({\bf x})\big)
d{\bf x}\right|.
\end{eqnarray}
Using condition (F.2) ($f(\cdot)$ is compactly supported),
$f(\cdot)$ is bounded away from zero on its support, thus, we have
for $n$ enough large,
there exists $\gamma>0$, such that
 $f({\bf x})>\gamma$, for all
${\bf x}$ in the support of $f(\cdot)$. By the same previous
arguments we have, for $n$ enough large,
\begin{eqnarray*}
\pi_n(\mathbb{J})&\leq&
\frac{1}{\gamma} \sup_{{\bf x}\in \mathbb{J}}
\left\vert  f_{n,h_n}({\bf x})- f({\bf x})
\right\vert\\&&+\Big(\frac{1}{\gamma^{2} }+1\Big)\sup_{{\bf x}\in
\mathbb{J}} \left\vert  f_{n,h_n}({\bf x})-
f({\bf x})\right\vert.
\end{eqnarray*}
One finds, by combining the last equation with
(\ref{Chap3-resulta-em-70}),
%e30 ###
\begin{eqnarray}\label{Chap3-resulta-em-701}
\lefteqn{\limsup_{n\rightarrow\infty} \sup_{h_n^\prime\leq h\leq h_n^{\prime\prime}}\frac{\sqrt{nh}~~\pi_n(\mathbb{J})}{
\sqrt{\{(\log(1/h)\vee
\log\log n)}}}\nonumber\\[-1pt]
&\leq& \left(\frac{\gamma^{2}+\gamma+1}{\gamma^{2}}\right)^{1/2}\sup_{{\bf x}\in \mathbb{J}}\left\{f({\bf x})\int_{\mathbb{R}^d}K^2({\bf t})d{\bf t}\right\}^{1/2}.
\end{eqnarray}
Let $\{\mathbb{J}_\ell\}$, $\ell=1,2,\ldots,$ be a sequence of
nondecreasing nonempty compact subsets of $\mathring{\mathbb{I}}$
such that
\vspace*{-3pt}
\begin{eqnarray*}
\bigcup_{\ell \geq 1} \mathbb{J}_\ell=\mathbb{I}.
\end{eqnarray*}
Now, from (\ref{Chap3-resulta-em-701}), it is straightforward to
observe that
\begin{eqnarray*}
\lefteqn{\lim_{\ell\rightarrow\infty}\limsup_{n\rightarrow\infty} \sup_{h_n^\prime\leq h\leq h_n^{\prime\prime}}\frac{\sqrt{nh\gamma_{n}^{4}}\pi_n(\mathbb{J}_\ell)}{
\sqrt{(\log(1/h)\vee
\log\log n)}}}\\[-1pt]
&\leq&
\lim_{\ell\rightarrow\infty}\left(\frac{\gamma^{2}+\gamma+1}{\gamma^{2}}\right)^{1/2}\sup_{{\bf x}\in \mathbb{J}_\ell}\left\{f({\bf x})\int_{\mathbb{R}^d}K^2({\bf t})d{\bf t}\right\}^{1/2}\\[-1pt]
&\leq&\left(\frac{\gamma^{2}+\gamma+1}{\gamma^{2}}\right)^{1/2}\sup_{{\bf x}\in \mathbb{I}}\left\{f({\bf x})\int_{\mathbb{R}^d}K^2({\bf t})d{\bf t}\right\}^{1/2}.
\end{eqnarray*}
\noindent The proof of Theorem \ref{Chap3-corol2} is completed.%\hfill$\blacksquare$

\vspace*{-3pt}
\subsection*{Proof of Theorem \ref{Chap3-theo}.}
\vspace*{-3pt}

Let $\varphi_{n,h_n}(\mathbf{x}):=\mathbb{E}(f_{n,h_n}(\mathbf{x}))$. Recall that
\begin{eqnarray*}
H_{n,h_n}^{(2)}(f)-\widehat{\mathbb{E}}H_{n,h_n}^{(2)}(f)\nonumber
&\!\!=\!\!&-\frac{1}{n}\sum_{i=1}^{n}
\mathds{1}_{\Omega_{n,i}}\log(f_{n,h_n}({\bf X}_i))+\mathds{1}_{\Omega_{n,i}}\log\left(\varphi_{n,h_n}({\bf
X}_{i})\right)\\[-1pt]
&\!\!=:\!\!&\boldsymbol{\Xi}_{n,h_n}.
\end{eqnarray*}
Using a Taylor-Lagrange expansion of the $\log(\cdot)$ function,
we have, for some random sequence $\theta_n\in(0,1)$,
\begin{equation*}
\boldsymbol{\Xi}_{n,h_n}= \frac{1}{n}\sum_{i=1}^{n}
\mathds{1}_{\Omega_{n,i}}\left[\frac{f_{n,h_n}({\bf X}_i)-\varphi_{n,h_n}({\bf
X}_{i})}{(1-\theta_n)f_{n,h_n}({\bf X}_i)+\theta_n\varphi_{n,h_n}({\bf
X}_{i})}\right].
\end{equation*}
Recalling that $\Omega_{n,i}=\big\{
f_{n,h_n}({\bf X}_i)\ge\gamma_{n}\big\}$,
we readily obtain, with probability 1,
\begin{eqnarray*}
\vert\boldsymbol{\Xi}_{n,h_n}\vert& \le & \frac{1}{n\gamma_{n}}\sum_{i=1}^{n} \mathds{1}_{\Omega_{n,i}}\left\vert f_{n,h_n}({\bf X}_i)-\varphi_{n,h_n}({\bf
X}_{i})\right\vert\\[-1pt]
& \le &\frac{1}{\gamma_{n}}\sup_{{\bf x}\in \mathbb{I}}
\left\vert f_{n,h_n}({\bf x})-\varphi_{n,h_n}({\bf
x})\right\vert \\[-1pt]
&=&\frac{1}{\gamma_{n}}\sup_{{\bf x}\in \mathbb{I}}
\left\vert f_{n,h_n}({\bf x})-\mathbb{E}(f_{n,h_n}({\bf
x}))\right\vert.
\end{eqnarray*}
Combining the last inequality with (\ref{Chap3-resulta-em}), we readily obtain the desired result.
%\hfill$\blacksquare$

%\noindent
\subsection*{Proof of Theorem \ref{Chap-5-theorem1}.}
We have
\begin{eqnarray}
H_{n,h_n}^{(2)}(f)-H(f)\nonumber &=&\{H_{n,h_n}^{(2)}(f)-\widehat{\mathbb{E}}H_{n,h_n}^{(2)}(f)\}+\{\widehat{\mathbb{E}}H_{n,h_n}^{(2)}(f)-H(f)\}.
\end{eqnarray}
Since the first term in the right hand of the last equality is controlled in the preceding proof, it remains only to evaluate the second one. To simplify our exposition, we will decompose
$\widehat{\mathbb{E}}H_{n,h_n}^{(2)}(f)-H(f)$ into the sum of three
components, that is
%e31 ###
\begin{eqnarray}
\widehat{\mathbb{E}}H_{n,h_n}^{(2)}(f)-H(f)&=&
-\frac{1}{n}\sum_{i=1}^{n}
\mathds{1}_{\Omega_{n,i}}\log(\varphi_{n,h_n}({\bf X}_i))+\mathbb{E}\left(\log\left(f({\bf
X}_{i})\right)\right)\nonumber\\
&=& -\frac{1}{n}\sum_{i=1}^{n}\mathds{1}_{\Omega_{n,i}}
\left(\log(\varphi_{n,h_n}({\bf X}_i))-\log(f({\bf X}_i))\right)\nonumber\\
&&-\frac{1}{n}\sum_{i=1}^{n}\left(\mathds{1}_{\Omega_{n,i}}
\log(f({\bf X}_i))-\log(f({\bf X}_i))\right)\nonumber\\
&&-\frac{1}{n}\sum_{i=1}^{n}\left(\log(f({\bf X}_i))-\mathbb{E}\left(\log(f({\bf X}_i))\right)\right)\nonumber\\
&=:& -\boldsymbol{\nabla}_{1,n,h_n} -\boldsymbol{\nabla}_{2,n,h_n}-\boldsymbol{\nabla}_{3,n,h_n}. \label{Chap-5-delta-def}
\end{eqnarray}
\noindent In view of (\ref{Chap-5-delta-def}), we have
\begin{eqnarray*}
\boldsymbol{\nabla}_{1,n,h_n}&=&\frac{1}{n}\sum_{i=1}^{n}
\mathds{1}_{\Omega_{n,i}}\left(\log(\varphi_{n,h}({\bf X}_i))-\log(f({\bf X}_i))\right).\\
\end{eqnarray*}
Using again a Taylor-Lagrange expansion of the $\log(\cdot)$ function,
we have, for some random sequence $\theta_n\in(0,1)$,
\begin{equation*}
\boldsymbol{\nabla}_{1,n,h_n}= \frac{1}{n}\sum_{i=1}^{n}
\mathds{1}_{\Omega_{n,i}}\left[\frac{\varphi_{n,h_n}({\bf X}_i)-f({\bf X}_i)}{(1-\theta_n)\varphi_{n,h_n}({\bf X}_i)+\theta_nf({\bf X}_i)}\right].
\end{equation*}
By condition (F.2), there exists a constant $\eta_{\mathbb{I}}>0$, such that $f({\bf x})>\eta_{\mathbb{I}}$ for all ${\bf x}\in \mathbb{I}$. It follows that for $n$ enough large that,  $f({\bf x})>\gamma_{n}$  for all ${\bf x}\in \mathbb{I}$.
Recalling that $\Omega_{n,i}=\big\{
f_{n,h_n}({\bf X}_i)\ge\gamma_{n}\big\}$,
we readily obtain, with probability 1,
\begin{eqnarray*}
\vert\boldsymbol{\nabla}_{1,n,h_n}\vert& \le & \frac{1}{n\gamma_{n}}\sum_{i=1}^{n} \mathds{1}_{\Omega_{n,i}}\left\vert \varphi_{n,h_n}({\bf X}_i)-f({\bf X}_i)\right\vert\\
& \le & \frac{1}{\gamma_{n}}\sup_{{\bf x}\in \mathbb{I}}
\left\vert \varphi_{n,h_n}({\bf x})-f({\bf x})\right\vert.
\end{eqnarray*}
We mention that the bandwidth $h$ is to be chosen in such a way that the bias of $f_{n,h}({\bf x})$ may be neglected, in the sense that
%e32 ###
\begin{eqnarray}\label{Chap-5-Ineq-bias}
&&\lim_{n\rightarrow\infty}\sup_{h_n^{\prime}\leq
h\leq h_n^{\prime\prime}}\left\{\frac{nh}{2\log(1/h)}
 \right\}^{1/2}\sup_{{\bf x}\in\mathbb{I}}\big|\varphi_{n,h}({\bf x})-f({\bf x})
 \big|=0,
\end{eqnarray}
which is implied by (\ref{equationbiais}).
Thus,
%e33 ###
\begin{equation}\label{Chap-5-Delta111}
\limsup_{n\rightarrow\infty} \sup_{h_n^{\prime}\leq h\leq
h_n^{\prime\prime}} \frac{\sqrt{nh\gamma_{n}^{2}}\vert\boldsymbol{\nabla}_{1,n,h}\vert}{\sqrt{2\log(1/h)}}=0.
\end{equation}
We next evaluate the second term $\boldsymbol{\nabla}_{2,n,h_n}$ in the
right side of (\ref{Chap-5-delta-def}). We have from (\ref{equationbiais}) and (\ref{Chap3-resulta-em})
\begin{eqnarray*}
\sup_{h^\prime_n\leq h\leq h^{\prime\prime}_n}\sup_{{\bf x}\in\mathbb{I}}\big|f_{n,h}({\bf x})-f({\bf x})
 \big|=O\left(\sqrt{\frac{(\log (1/h^\prime_n)}{nh^\prime_n}}\right).
\end{eqnarray*}
Thus, for $n$ sufficiently large, almost surely, $f_{n,h}({\bf x})\geq (1/2)f({\bf x})$ for all ${\bf x}\in
\mathbb{I}$ and all $h\in[h_n^{\prime},h_n^{\prime\prime}]$.
Note that under condition (F.2), the density $f(\cdot)$ is compactly supported, it is
 possible to find a positive constant $\eta_\mathbb{I}$ such as $f(\mathbf{x})>\eta_\mathbb{I}$.
This implies that $f_{n,h}({\bf x})\geq \eta_\mathbb{I}/2$, and thus, for all $n$ enough large,
we have, almost surely,
%e34 ###
\begin{equation}
    \mathds{1}_{\Omega_{n,i}}=1,
\end{equation}
which implies that, for all $n$ enough large, almost surely,
%e35 ###
\begin{equation}\label{Chap-5-ineq13}
    \boldsymbol{\nabla}_{2,n,h_n}=0.
\end{equation}
\noindent We finally evaluate the second term $\boldsymbol{\nabla}_{3,n,h_n}$ in the
right side of (\ref{Chap-5-delta-def}). We have,
\begin{eqnarray*}
\boldsymbol{\nabla}_{3,n,h_n}&=&-\frac{1}{n}\sum_{i=1}^{n}\boldsymbol{\xi}_{i}\nonumber,
\end{eqnarray*}
where, for $i=1,\dots,n$,  $$\boldsymbol{\xi}_{i}:=\log\{f({\bf X}_i)\}
-\mathbb{E}\Big(\log\{f({\bf X}_i)\}\Big),$$
are a centered independent and identically distributed  random variables with finite variance
$\mbox{Var}\big(\log(f({\bf X}_i))\big)$ (condition (F.3)). Observe that
\begin{eqnarray*}\label{Chap-5-ineq12}
\nonumber\frac{\gamma_{n}}{n}\frac{\sqrt{nh_n}
\sum_{i=1}^{n} \boldsymbol{\xi}_{i}}{\sqrt{2\log(1/h_n)}}
&=&\frac{\gamma_{n}\sqrt{h_n\log\log
n}}{\sqrt{\log(1/h_n)}}
\frac{\sum_{i=1}^{n} \boldsymbol{\xi}_{i}}{\sqrt{2n\log\log n}}
\end{eqnarray*}
which, by the law of the iterated logarithm, tends to $0$ as
$n$ tends to infinity. Namely,
%e36 ###
\begin{eqnarray}\label{Chap-5-ineq120}
\lim_{n\rightarrow\infty} \sup_{h_n^{\prime}\leq h\leq
h_n^{\prime\prime}} \frac{\sqrt{nh\gamma_{n}^{2}}\vert\boldsymbol{\nabla}_{3,n,h}\vert}{\sqrt{2\log(1/h)}}=0.
\end{eqnarray}
Using (\ref{Chap-5-ineq120}) and (\ref{Chap-5-ineq13}) in connection with (\ref{Chap3-resulta-em-70})
completes the proof of Theorem \ref{Chap-5-theorem1}. %\hfill$\blacksquare$

%\section*{Acknowledgments}
%The authors are grateful to an anonymous referee and an Associate Editor for several
%insightful suggestions that  led to a significant improvement of the presentation and the correction
%of a technical argument in the original proof.

%\def\cprime{$'$}

\end{document}